\documentclass[10pt]{article}
\usepackage{mathrsfs}
\usepackage{stmaryrd}
\usepackage{amsfonts,amsmath,amssymb,amscd}
\usepackage[colorlinks,linkcolor=red,anchorcolor=blue,citecolor=green]{hyperref}
\usepackage{shadow}
\usepackage{graphics,graphicx}
\usepackage{color}
\usepackage{enumerate}
\usepackage{longtable}
\usepackage{appendix}
\allowdisplaybreaks

\newtheorem{thm}{Theorem}[section]
\newtheorem{lem}[thm]{Lemma}

\newtheorem{defi}[thm]{Definition}
\newtheorem{remark}[thm]{Remark}

\newtheorem{conj}[thm]{Conjecture}

\def\qed{\hfill \rule{4pt}{7pt}}

\def\pf{\noindent {\it{Proof.}\hskip 2pt}}

\begin{document}
\begin{center}
{\large\bf Cutpoints for Random Walks on Quasi-Transitive Graphs}
\end{center}

\begin{center}
He Song and Kai-Nan Xiang

\footnotesize{School of Mathematical Sciences, LPMC, Nankai University}\\
\footnotesize{Tianjin City, 300071, P. R. China}\\
\footnotesize{Emails: songhe@mail.nankai.edu.cn} (H. Song)\\
\footnotesize{~~~~~~~~~~~~kainanxiang@nankai.edu.cn} (K. N. Xiang)
\footnote{The project is supported partially by CNNSF (No. 11271204).

   {\it MSC2010 subject classifications}. 60J10, 05C81, 60D05, 60G17.

   }
\end{center}
\begin{abstract}
We prove that a simple random walk on quasi-transitive graphs with the volume growth
being faster than any polynomial of degree 4 has a.s. infinitely many cut times, and hence infinitely many cutpoints.
This confirms a conjecture raised by I. Benjamini,
 O. Gurel-Gurevich and O. Schramm [2011, Cutpoints and resistance of random walk paths, {\it Ann. Probab.} {\bf 39(3)}, 1122-1136]
 that PATH of
 simple random walk on any transient vertex-transitive graph has a.s. infinitely many cutpoints in the corresponding case.
\end{abstract}

{\large\bf Claim: Chinese version has been published in Acta Mathematica Sinica, Chinese Ser. 60 (2017), no.6, 947-954.}

\section{Introduction}
Let $G=(V, E)$ be a locally finite, connected infinite graph, and $S=\{S_n\}_{n=0}^\infty$ a random walk (RW) on $G$
with each $S_{n+1}$ being a neighbor of $S_n.$ Define path graph {\it PATH} of $S$ as the subgraph consisting of all vertices and all edges visited by it.

\begin{defi}
For a RW $S=\{S_n\}_{n=0}^\infty$ on $G$, a nonnegative integer $n$ is called a cut time for $S$ if
\[S[0,n]\cap S[n+1,\infty )=\emptyset,\]
where $S[0,n]=\{S_j:0\leq j\leq n\}$ and $S[n+1,\infty)=\{S_j:j\geq n+1\}$.

A vertex in PATH is called a cutpoint of PATH if it separates $S_0$ from infinity; namely, if we delete it from PATH, then it results $S_0$ in a finite
connected component.
\end{defi}

Note if nonnegative integer $n$ is a cut time, then the vertex $S_n$ must be a cutpoint of PATH; while the converse is not true.
And recurrent RWs have no cut time. So the cut time is only interesting for transient RWs.
Transience doesn't imply the infinity of the number of cut times for RWs \cite{NJRL2008}; whereas
I. Benjamini, O. Gurel-Gurevich and O. Schramm \cite{IBOGG2011} proved that for every transient Markov chain, the expected number of cut times
is infinite. Recall that for any transient RW on $G,$ PATH is always a recurrent graph (see \cite{IBOGG2005}, \cite{IBOGG2007} and also \cite{RLYP2013} Section 9.5). A natural and interesting question is that on what kind of graphs the PATH has infinitely many cutpoints.

To continue, let $Aut(G)$ be the set of all automorphisms of $G$; and notice that $G$ is called vertex transitive if $Aut(G)$ acts transitively on $G$, and quasi-transitive if $Aut(G)$ acts with finitely many orbits. In \cite{IBOGG2011}, I. Benjamini,
 O. Gurel-Gurevich and O. Schramm raised the following conjecture:

\begin{conj}\label{conjecture1}
The PATH of simple random walk (SRW) $S=\{S_n\}_{n=0}^{\infty}$ on any transient vertex transitive graph $G$ has a.s. infinitely many cutpoints.
\end{conj}

For any vertex $x$ of $G,$ let $B(x,n)$ be the ball in $G$ with radius $n$ centered at $x,$ and $\vert B(x,n)\vert$ the cardinality of $B(x,n).$
Fix a vertex $v_0$ of $G,$ let
$$V(n)=V_G(n)=\vert B(v_0,n)\vert.$$
In this paper, we prove Conjecture \ref{conjecture1} in the case of (1.1) by proving

\begin{thm}\label{Thm2}
Given any  quasi-transitive infinite graph $G$ with the volume growth satisfying
 $$\limsup_{n\rightarrow\infty} \frac{V(n)}{n^4}=\infty. \eqno{(1.1)}$$
Then almost surely there are infinitely many cut times, and hence infinitely many cutpoints, for SRW on $G$.
\end{thm}

Note the study of cut times for RWs on $\mathbb{Z}^d$ can date back to P. Erd\"{o}s and S. J. Taylor \cite{PESJT1960}.
G. Lawler \cite{GL1991} studied cut times of SRW on $\mathbb{Z}^d$ and proved that
SRW has a.s. infinitely many cut times when $d\geq 4.$ Then N. James and Y. Peres \cite{NJYP1996} solved the case of $d=3$. Thus every transient SRW on $\mathbb{Z}^d$ has a.s. an infinite number of cut times. Indeed \cite{NJYP1996} proved that the cut time number of every transient RW on $\mathbb{Z}^d$ with finite range is infinite almost surely. Furthermore, G. Lawler \cite{GFL1996} verified that a bilateral random walk in $\mathbb{Z}^d$ has almost surely infinitely many cut times if
$d\geq 5;$ whereas if $d\leq 4,$ there are almost surely none. For a criterion for a.s. finiteness and infiniteness of the number of cutpoints (cut times)
in case of transient RWs on $\mathbb{Z}_{+},$ see E. Cs\'{a}ki, A. F\"{o}ldes and P. R\'{e}v\'{e}sz \cite{ECAF2010}.

Turn to finitely generated groups $G$ with polynomial volume growth of degree $D.$ When $D\leq 2,$ every symmetric RW on $G$ is recurrent.
When $D\geq 5$, N. James and Y. Peres  \cite{NJYP1996} proved that a.s. there are infinite number of cut times for symmetric irreducible RW with finite range;
and G. Alexopoulos \cite{GA2002} proved that this is still true without irreducibility. Note for any $D,$ for any non-symmetric RW on $G,$
the number of cut times is a.s. infinite (\cite{GA2002}). When $D=3,4,$ for symmetric RWs on $G$, S. Blach\`{e}re \cite{BS2003} verified that there are a.s.
infinitely many cut times.

\begin{remark}\label{remark1}
For a quasi-transitive infinite graph $G,$ the followings hold (see Appendix \ref{appendix}).

{\bf (i)} If $G$ has polynomial volume growth, then there must be a natural number $D$ satisfying
\[0<\liminf_{n\rightarrow\infty}\frac{V(n)}{n^D}\leq \limsup_{n\rightarrow\infty}\frac{V(n)}{n^D}<\infty;\]
which implies that (1.1) is equivalent to
$\limsup\limits_{n\rightarrow\infty} \frac{V(n)}{n^5}>0.$

{\bf (ii)} SRW on $G$ is transient if and only if $\displaystyle \limsup_{n\rightarrow\infty} \frac{V(n)}{n^3}>0.$

{\bf (iii)} Following the same line as the proof for finitely generated groups $G$  by N. James and Y. Peres  \cite{NJYP1996}, the similar results can be extended to transitive graphs,
\end{remark}
\vskip 3mm

The left case for Conjecture \ref{conjecture1} is that  $G$ has polynomial volume growth with degree $D=3$ or $4$.
Though our proof does not work for the case, we believe SRW on $G$ still has a.s. infinitely many cut times in the case.
In addition, for any reversible RW $S=\{S_n\}_{n=0}^{\infty}$ on $G$ such that for any $x,y\in V,$
\begin{eqnarray*}
p(x,y)=p(\gamma (x),\gamma (y)),\ \forall\gamma\in Aut(G),\ \mbox{and}\ p(x,y)>0\Longleftrightarrow y \ \mbox{is adjacent to}\ x,
\end{eqnarray*}
where $p(x,y)$ is 1-step transition probability for $S$ from $x$ to $y,$
Theorem 1.3 is true by following the same line as
the proof in this paper.\\

\section{Proof of Theorem \ref{Thm2}}

To begin our proof, let us firstly introduce some preliminaries.
Let $S^1=\{S^1_n\}_{n=0}^{\infty}$ and $S^2=\{S^2_n\}_{n=0}^{\infty}$ be two independent SRWs starting at a given vertex $o$ in $G$. Extend SRW $S^1$ to a two-sided SRW $S=\{S_n\}_{n\in\mathbb{Z}}$:
\begin{align*}
S_j=
 \left\{
   \begin{array}{ll}
     S^1_j, & 0\leq j<\infty, \\
     S^2_{-j}, & -\infty<j< 0.
   \end{array}
 \right.
\end{align*}
\vskip 3mm

\begin{defi}
Call a time $j$ loop-free, if
$$S(-\infty,j]\cap S[j+1,\infty)=\emptyset.$$
Here $$S(-\infty,j]=\{S_n:\ n\leq j\},\ S[j+1,\infty)=\{S_n:\ n\geq j+1\}.$$
\end{defi}
\vskip 3mm

Clearly, if a time $j\geq 0$ is loop-free, then $j$ is a cut time for $S^1.$\\

\begin{defi}
Given any $i,j\geq 0.$  Call $(i,j)$ a $*$-last intersection, if
 $$S^1_i=S^2_j,\ S^1_{i_1}\neq S^2_{j_1},\ (i,j)<(i_1,j_1)<(\infty,\infty).$$
Here $(i_1,i_2)\leq (j_1,j_2)$ means $i_1\leq j_1$ and $i_2\leq j_2$; and $(i_1,i_2)<(j_1,j_2)$ means $(i_1,i_2)\leq (j_1,j_2)$ but $(i_1,i_2)\neq (j_1,j_2).$
\end{defi}
\vskip 3mm

Note that the $*$-last intersection may not be unique.


%
%

\begin{lem}\label{lem6} If quasi-transitive $G$ satisfies that for some positive constant $C$ and natural number $D,$
 $$V(n)\geq Cn^D,\ \forall n\geq 1;$$
then
\[\sup_{x,y\in V}\frac{p^{(n)}(x,y)}{\deg (y)} \leq n^{-\frac{D}{2}},\]
where $p^{(n)}(x,y)$ is $n$-step transition probability for the SRW from $x$ to $y.$
\end{lem}
\vskip 3mm

\begin{lem}\label{lem7}
For quasi-transitive $G$ satisfying (1.1), $\displaystyle \sum_{j=0}^\infty jp^{(j)}(o,o)<\infty$.
\end{lem}
\pf\
By Remark 1.4(i) and (1.1), we have that $V(n)\geq C n^5$ for some positive constant $C.$
Therefore, by Lemma \ref{lem6}, there exists a positive constant $c_1$ such that
\[p^{(n)}(x,y)\leq c_1n^{-\frac{5}{2}},\ x,y\in V,\ n\geq 1.\]
And further,
$$\displaystyle \sum_{j=0}^\infty jp^{(j)}(o,o)\leq c_1\sum\limits_{j=0}^{\infty}\frac{1}{n^{3/2}}<\infty.$$
\qed
\vskip 4mm

Let
\[g(o)=\mathbb{P}\left[S^1_i\neq S^2_j,\ (0,0)<(i,j)<(\infty,\infty)\right]. \eqno{(2.1)}\]
Suppose quasi-transitive graph $G=(V,E)$ has $k$ types of vertices. That is we can divide the vertex set $V$ into $k$ orbits according to
$Aut(G)$-action on $V:$ $[1]$, $[2]$, $\cdots$, $[k]$.
Clearly, if initial points $o_1$ and $o_2$ are in the same orbit $[i],$ then
$$g(o_1)=g(o_2).$$
Thus function $g(\cdot)$ takes at most $k$ values.
Let $\widetilde{g}=\displaystyle\max_{o\in V}g(o)$ and $R$ be the number of intersection times, namely
\[R=\displaystyle \sum_{i=0}^\infty \sum_{j=0}^\infty I_{\left\{S^1_i=S^2_j\right\}}.\]
Notice that
\[\mathbb{E}(R)=\displaystyle \sum_{i=0}^\infty \sum_{j=0}^\infty \mathbb{P}\left\{S^1_i=S^2_j\right\}.\]\\

\begin{lem}\label{Lem4}
For quasi-transitive $G=(V,E)$ satisfying (1.1),
\[\widetilde{g}\geq c:=1/\mathbb{E}(R)>0.\]
\end{lem}
\vskip 3mm
\pf\  Since for any $x,y\in V$ and any $n\geq 1,$
$$\deg (x)\ p^{(n)}(x,y)=\deg (y)\ p^{(n)}(y,x),$$
we have that
\begin{align*}
\mathbb{P}(S^1_i=S^2_j)&=\displaystyle \sum_{x\in V}\mathbb{P}(S^1_i=x,~S^2_j=x)=\displaystyle \sum_{x\in V}\mathbb{P}(S^1_i=x)\mathbb{P}(S^2_j=x)\\
&=\displaystyle \sum_{x\in V}p^{(i)}(o,x)p^{(j)}(o,x)=\displaystyle \sum_{x\in V}p^{(i)}(o,x)\frac{\deg(x)}{\deg(o)}p^{(j)}(x,o).
\end{align*}
Note $d=\max\limits_{v\in V}\deg (v)<\infty.$ Then
\[\mathbb{P}(S^1_i=S^2_j)\leq \displaystyle d\sum_{x\in V} p^{(i)}(o,x)p^{(j)}(x,o)=d\ p^{(i+j)}(o,o).\]
Therefore,
\[\mathbb{E}(R)\leq \displaystyle d\sum_{i=0}^\infty \sum_{j=0}^\infty p^{(i+j)}(o,o)=\displaystyle d\sum_{\ell=0}^\infty (\ell
+1)p^{(\ell)}(o,o).\]
By Lemma \ref{lem7}, we see that
$$\mathbb{E}(R)<\infty.$$
Thus, a.s. the number of intersections for $S^1$ and $S^2$ is finite. This implies the existence of a $*$-last intersection for $S^1$ and $S^2$
almost surely, and
\[1\leq\displaystyle \sum_{i=0}^\infty \sum_{j=0}^\infty\mathbb{P}((i,j)~ is~ a~ *-last~ intersection).\]
But
\begin{align*}
&\mathbb{P}((i,j)~ is~ a~ *-last~ intersection)\\&=\mathbb{P}\left\{S^1_i=S^2_j,~ S^1_{i_1}\neq S^2_{j_1},\ (i,j)<(i_1,j_1)< (\infty,\infty)\right\}\\
&=\mathbb{P}\left\{S^1_i=S^2_j\right\}\mathbb{P}\left\{\left. S^1_{i_1}\neq S^2_{j_1},\ (i,j)<(i_1,j_1)< (\infty,\infty)\right\vert S^1_i=S^2_j\right\}.
\end{align*}
Notice (2.1). By the Markov property for $S^1$ and $S^2,$ we obtain that
\[\mathbb{P}\left\{\left.S^1_{i_1}\neq S^2_{j_1},\ (i,j)<(i_1,j_1)< (\infty,\infty)\right\vert S^1_i=S^2_j\right\}=g\left(S^1_i\right)
     =g\left(S^2_j\right).\]
Thus,
\[1\leq\displaystyle \sum_{i=0}^\infty \sum_{j=0}^\infty\mathbb{P}((i,j)~ is~ a~ *-last~ intersection)\leq \sum_{i=0}^\infty \sum_{j=0}^\infty\mathbb{P}\left(S^1_{i}=S^2_{j}\right)\widetilde{g}=\widetilde{g}\ \mathbb{E}(R).\]
Namely, $\widetilde{g}\geq 1/\mathbb{E}(R)>0.$ \qed
\vskip 4mm

Without loss of generality, assume that
$$[1]\subseteq \left\{v\in G\left|~g(v)=\widetilde{g}\right.\right\}.  \eqno{(2.2)}$$
Since $G$ is infinite and has $k$ orbits, at least one orbit has infinitely many vertices. In fact, we have\\

\begin{lem}
Every orbit has infinitely many vertices.
\end{lem}

\pf \ Suppose that there exists an orbit, say orbit $[k],$ with finitely many vertices. As said above, there is an orbit, say $[i],$
having infinitely many vertices. Choose a vertex $x\in [i]$ and a large enough natural number $r$ such that $[k]\subset B(x,r).$ Choose another vertex $y\in [i]$ such that the distance between $x$ and $y$ on $G$ is bigger than $2(r+1).$ Then $B(y,r)\cap [k]=\emptyset.$ But $x$ and $y$ are in the same orbit $[i],$
it is impossible that
   $$[k]\subset B(x,r),\ B(y,r)\cap [k]=\emptyset.$$
A contradiction! Therefore, every orbit has infinitely many vertices.\qed
\vskip 4mm

\begin{lem}\label{Lem8}
The SRW on $G$ visits every orbit infinite times.
\end{lem}

\pf\
Consider the SRW $S^1=\{S_n^1\}_{n=0}^{\infty}$ on $G$. Define a random walk $X=\{X_n\}_{n=0}^{\infty}$ on $\{1,2,\cdots, k\}$ as follows:
$$X_n=i\ \mbox{if}\ S_n^1\in [i],\ n\geq 0.$$
For any $n\geq 0,$ and any $i,j\in\{1,2,\cdots,k\},$ and any $x\in [i],$  it is easy to see that
$$\mathbb{P}\left[\left. S_{n+1}^1\in [j]\ \right\vert S^1_{n}=x\right]\ \mbox{does not depend on}\ x\ \mbox{and}\ n,\  \mbox{and is a function in}\ i,j.$$
Let $p_{i,j}=\mathbb{P}\left[\left. S_{n+1}^1\in [j]\ \right\vert S^1_{n}=x\right].$ Then we get that $X$ is a Markov chain on $\{1,2,\cdots,k\}$ with transition matrix $(p_{i,j})_{1\leq i,j\leq k}.$ Due to $G$ is connected, for any $1\leq i,j\leq k,$ we see that
there is a natural number $n$ satisfying
$$\mathbb{P}\left[\left. S_n^1\in [j]\ \right\vert S_0^1\in [i]\right] >0;$$
which means $X$ is irreducible. Notice that $\{1,2,\cdots,k\}$ is finite and $X$ is irreducible. By the standard theory of Markov chains,
we obtain immediately that every state $i\in\{1,2,\cdots,k\}$ is recurrent for $X.$
Thus SRW on $G$ visits every orbit infinite times.\qed\\

Note (2.1)-(2.2). Then Lemma \ref{Lem4} implies that for the two-sided SRW $S$ with $S(0)=o\in [1],$
\[\mathbb{P}(S(0,\infty)\cap S(-\infty,0]=\emptyset)=\mathbb{P}\left(S^1_i\neq S^2_j,~(0,0)<(i,j)<(\infty, \infty)\right)=c>0.  \eqno{(2.3)}\]
For arbitrary $o\in V,$ by Lemma \ref{Lem8}, we have that the following stopping times are all finite almost surely:
\[\tau_0=\inf\left\{r\geq 0:\ S^1_r\in [1]\right\},\ \tau_{n+1}=\inf\left\{r>\tau_n:\ S^1_r\in [1]\right\},\ n=0,1,2,\cdots.\]
Clearly, almost surely, as $n\uparrow\infty,$ $\tau_n\uparrow \infty.$ Let
\[cut_n=\{\tau_n~ is~a~ cut~ time\},\ n=0,1,2,\cdots.\]
Then for any $n\geq 0,$
\[\mathbb{P}(cut_n)=\mathbb{P}\left(\left\{S^1_0,S^1_1,\cdots, S^1_{\tau_n}\right\}\cap\left\{S^1_{\tau_n+1},\cdots,S^1_\infty\right\}=\emptyset\right).
    \eqno{(2.4)}\]

\vskip 4mm
\begin{lem}\label{Lem2}
Suppose $o\in [1].$ Then with probability one, the SRW $S^1=\left\{S^1_n\right\}_{n=0}^{\infty}\ \left(S_0^1=o\right)$ on quasi-transitive $G$ satisfying (1.1) has infinitely many cut times $m$ with $S_m^1\in [1].$
\end{lem}

\pf\ Fix an arbitrary natural number $n.$ For any $y\in [1],$ there is a $\Phi_{y}\in Aut(G)$ such that $\Phi_{y}(y)=o.$ Notice that
the image $\widetilde{S}^1=\left\{\widetilde{S}_0^1,\widetilde{S}_1^1,\cdots\right\}$ of path $\left\{S^1_{\tau_n},S^1_{\tau_n+1},\cdots\right\}$ under $\Phi_{S_{\tau_n}^1}$ is a trajectory of an SRW
with initial point $o.$
Consider the following path $\widetilde{S}^2:=\left\{\widetilde{S}_0^2,\widetilde{S}_1^2,\cdots, \widetilde{S}_{\tau_n}^2\right\}:$
$$\widetilde{S}_{\ell}^2=\Phi_{S_{\tau_n}^1}(S_{\tau_n-\ell}^1),\ \ell=0,1,2,\cdots, \tau_n.$$
For any finite path $p=(w_0,w_1,\cdots,w_{\ell})$ in $G,$ let
$\widetilde{p}=(w_\ell,w_{\ell-1},\cdots,w_1,w_0),$ and
$$f(p)=\prod\limits_{j=0}^{\ell-1}\frac{1}{\deg (w_j)},\
  \widetilde{f}(p)=f\left(\widetilde{p}\right)=\prod\limits_{j=0}^{\ell-1}\frac{1}{\deg (w_{\ell-j})}.$$
Let $\widetilde{P}(n)$ be the trajectory space for $\left\{\widetilde{S}_0^2,\widetilde{S}_1^2,\cdots, \widetilde{S}_{\tau_n}^2\right\}.$
For any $x\in V$ with $Q_x:=\mathbb{P}\left[S^1_{\tau_n}=x\right]>0,$
write $P(n,x)$ for the corresponding trajectory space for $\left(S^1_0, S^1_1,\cdots, S^1_{\tau_n}\right)$ with $S_{\tau_n}^1=x.$
Then for any $p\in \widetilde{P}(n),$ given $S_{\tau_n}^1=x\in [1]$ with $Q_x>0,$
the conditional probability of
$$\left\{\left\{S_0^1,S_1^1,\cdots,S^1_{\tau_n}\right\}=\Phi_{S_{\tau_n}^1}^{-1}\circ \widetilde{p}\right\}$$
is
$\frac{f\left(\widetilde{p}\right)}{Q_x}I_{\left\{\Phi_{x}^{-1}\circ \widetilde{p}\in P(n,x)\right\}}.$
Note that both initial point and end point for path $p$ are in $[1],$ and thus have the same degree, $f\left(\widetilde{p}\right)=f(p).$
Therefore, for any $p\in \widetilde{P}(n),$
\begin{eqnarray*}
&&\mathbb{P}\left[\left\{\widetilde{S}_0^2,\widetilde{S}_1^2,\cdots, \widetilde{S}_{\tau_n}^2\right\}=p\right]
  =\mathbb{P}\left[\{S_0^1,S_1^1,\cdots,S^1_{\tau_n}\}=\Phi_{S_{\tau_n}^1}^{-1}\circ \widetilde{p}\right]\\
&&\ \ \ \ =\mathbb{E}\left[\mathbb{E}\left[\left. I_{\left\{\left\{S_0^1,S_1^1,\cdots,S^1_{\tau_n}\right\}=\Phi_{S_{\tau_n}^1}^{-1}\circ \widetilde{p}\right\}}\
                \right\vert S_{\tau_n}^1\right]\right]\\
&&\ \ \ \ =\sum\limits_{x}Q_x\ \frac{f\left(\widetilde{p}\right)}{Q_x}I_{\left\{\Phi_{x}^{-1}\circ \widetilde{p}\in P(n,x)\right\}}\\
&&\ \ \ \ =\sum\limits_xf\left(\widetilde{p}\right)I_{\left\{\Phi_{x}^{-1}\circ \widetilde{p}\in P(n,x)\right\}}=f\left(\widetilde{p}\right).
\end{eqnarray*}
Therefore,
\[\mathbb{P}\left[\left\{\widetilde{S}_0^2,\widetilde{S}_1^2,\cdots, \widetilde{S}_{\tau_n}^2\right\}=p\right]=f\left(\widetilde{p}\right)=f(p),\]
and $\widetilde{S}^2=\left\{\widetilde{S}_0,\widetilde{S}_1,\cdots, \widetilde{S}_{\tau_n}\right\}\in \widetilde{P}(n)$
evolves as an SRW path starting at $o.$

By the strong Markov property for the SRW, given $S_{\tau_n}^1,$ $\left\{S_{\tau_n+m}^1\right\}_{m\geq 0}$ and $\left\{S^1_m\right\}_{0\leq m\leq \tau_n}$
are independent. Hence, $\widetilde{S}^1$ and $\widetilde{S}^2$ are independent. Combining this with that $\widetilde{S}^1$ is an SRW with initial point $o,$
and $\widetilde{S}^2$ evolves as an SRW starting at $o;$ by (2.3) and (2.4), we have
\begin{eqnarray*}
\ \ \ \ \mathbb{P}(cut_n)&=&\mathbb{P}\left(\left\{S^1_0,S^1_1,\cdots, S^1_{\tau_n}\right\}\cap\left\{S^1_{\tau_n+1},\cdots,S^1_\infty\right\}=\emptyset\right)
                    =\mathbb{P}\left(\widetilde{S}^2\cap \left\{\widetilde{S}^1_1,\widetilde{S}^1_2,\cdots\right\}=\emptyset\right)\\
               &\geq &\mathbb{P}(S(0,\infty)\cap S(-\infty,0]=\emptyset)=\mathbb{P}\left(S^1_i\neq S^2_j,~(0,0)<(i,j)<(\infty, \infty)\right)=c>0.
               \ \ \ \  (2.5)
\end{eqnarray*}
Clearly, the above inequality also holds for $n=0.$

Thus by the Kochen-Stone lemma \cite{SKCS1964},
\[\mathbb{P}(cut_n~ i.o.)\geq \limsup_{m\rightarrow \infty} \frac{\sum_{0\leq i,j\leq m}\mathbb{P}(cut_i)\mathbb{P}(cut_j)}{\sum_{0\leq i,j\leq
  m}\mathbb{P}(cut_i\cap cut_j)}\geq c^2.\eqno{(2.6)}\]
Here $i.o.$ stands for infinitely often, and
we have used the following obvious inequalities:
$$\frac{\mathbb{P}(cut_i)\mathbb{P}(cut_j)}{\mathbb{P}(cut_i\cap cut_j)}\geq c^2,\ \forall 0\leq i,j<\infty.$$

Since $S^1=\left\{S_m^1\right\}_{m\geq 0}$ is transient, given path $p_{\tau_n}=\left\{S^1_0,S^1_1,\cdots,S^1_{\tau_n}\right\}$ with $n\geq 1,$
when natural number $n_0$ is large enough, $\left\{S^1_{n_0+m}\right\}_{m\geq 0}$ avoids $p_{\tau_n}.$ Therefore,
\begin{eqnarray*}
&&\limsup\limits_{m\rightarrow\infty}\mathbb{P}\left(\left.\left\{\tau_{n+m}\ is\ a\ cut\ time\ for\ \left\{S^1_{\tau_n+\ell}\right\}_{\ell\geq 0},\
      \left\{S^1_{\tau_{n+m}+\ell}\right\}_{\ell\geq 1}\ intersects\ p_{\tau_n}\right\}\right\vert p_{\tau_n}\right)\\
&&\ \ \ \ \ \ \ \ \ \ \leq \limsup\limits_{m\rightarrow\infty}\mathbb{P}\left(\left.\left\{\left\{S^1_{\tau_{n+m}+\ell}\right\}_{\ell\geq 1}\ intersects\
                            p_{\tau_n}\right\}\right\vert p_{\tau_n}\right)\\
&&\ \ \ \ \ \ \ \ \ \ \leq \mathbb{P}\left(\left.\limsup\limits_{m\rightarrow\infty}\left\{\left\{S^1_{\tau_{n+m}+\ell}\right\}_{\ell\geq 1}\ intersects\
                            p_{\tau_n}\right\}\right\vert p_{\tau_n}\right)\\
&&\ \ \ \ \ \ \ \ \ \ =0.
\end{eqnarray*}
Notice that given $p_{\tau_n},$ $\left\{S^1_{\tau_n+\ell}\right\}_{\ell\geq 0}$ is an SRW starting from $S^1_{\tau_n}\in [1];$
and (2.5) holds for any $o\in [1].$ Thus for any $m\geq 1,$
$$\mathbb{P}\left(\left.\left\{\tau_{n+m}\ is \ a \ cut\ time\ for\  \left\{S^1_{\tau_n+\ell}\right\}_{\ell\geq 0}\right\}\right\vert p_{\tau_n}\right)\geq c.$$
And further,
\begin{eqnarray*}
&&\liminf\limits_{m\rightarrow\infty}\mathbb{P}(cut_{n+m}\ \vert\ p_{\tau_n})\\
&&\geq \liminf\limits_{m\rightarrow\infty}\mathbb{P}\left(\left.\left\{\tau_{n+m}\ is\ a\ cut\ time\ for\
                           \left\{S^1_{\tau_n+\ell}\right\}_{\ell\geq 0}\right\}\right\vert p_{\tau_n}\right)\\
&&\ \ \ \ -\limsup\limits_{m\rightarrow\infty}\mathbb{P}\left(\left.\left\{\tau_{n+m}\ is\ a\ cut\ time\ for\
                              \left\{S^1_{\tau_n+\ell}\right\}_{\ell\geq 0},\ \left\{S^1_{\tau_{n+m}+\ell}\right\}_{\ell\geq 1}\ intersects\ p_{\tau_n}\right\}\right\vert p_{\tau_n}\right)\\
&&=\liminf\limits_{m\rightarrow\infty}\mathbb{P}\left(\left.\left\{\tau_{n+m}\ is\ a\ cut\ time\ for\
                           \left\{S^1_{\tau_n+\ell}\right\}_{\ell\geq 0}\right\}\right\vert p_{\tau_n}\right)\\
&&\geq c.
\end{eqnarray*}
Choose natural number $m_0$ large enough such that
$$\mathbb{P}(cut_{n+m}\ \vert\ p_{\tau_n})\geq c/2, \ \forall m\geq m_0.$$
Similarly to (2.6), by the Kochen-Stone lemma \cite{SKCS1964},
$$\mathbb{P}(cut_m\ i.o.\ \vert\ p_{\tau_n})\geq \mathbb{P}(cut_{n+m_0+m}\ i.o.\ \vert\ p_{\tau_n})\geq \frac{c^2}{4}.$$
Namely,
\[\mathbb{P}\left(cut_m~ i.o.\ \left|\ S^1_0, S^1_1,\cdots, S^1_{\tau_n}\right.\right)\geq \frac{c^2}{4}.\eqno{(2.7)}\]

By the L\'{e}vy 0-1 law, almost surely, as $n\rightarrow\infty,$
\[\mathbb{E}\left(I_{\{cut_m~ i.o.\}}\ \left|\ S^1_0, S^1_1,\cdots, S^1_{\tau_n}\right.\right)\rightarrow I_{\{cut_m~ i.o.\}}.\]
Thus
\[\mathbb{P}(cut_m~ i.o.)=1.\]
Namely the SRW $S^1$ has a.s. infinitely many cut times $m$ with $S_m^1\in [1].$
\qed\\

\vskip 4mm
\begin{lem}\label{Lem3}
For any $o\not\in [1],$ with probability one, the SRW $S^1=\left\{S^1_n\right\}_{n=0}^{\infty}\ \left(S_0^1=o\right)$ on quasi-transitive $G$ satisfying (1.1) has infinitely many cut times.
\end{lem}

\pf\  Notice almost surely $\tau_0\in (0,\infty).$ Let
$\widehat{S}^1=\left(S^1_{\tau_0+m}\right)_{m\geq 0}.$
Then $\widehat{S}^1$ is an SRW starting at $S^1_{\tau_0}\in [1].$ By Lemma \ref{Lem2}, almost surely,
there are infinitely many cut times, denoted in increasing order by $\sigma_1,\sigma_2,\cdots,$ for $\widehat{S}^1,$
such that $\tau_0+\{\sigma_1,\sigma_2,\cdots\}\subseteq\{\tau_0,\tau_1,\tau_2,\cdots\}.$
Note for any $n\geq 1,$
$$\left\{S_{\tau_0+\sigma_n+1}^1,S_{\tau_0+\sigma_n+2}^1,\cdots\right\}\cap \left\{S_{\tau_0}^1,S_{\tau_0+1}^1,\cdots,S^1_{\tau_0+\sigma_n}\right\}=\emptyset.$$
Let $A_n=\left\{\tau_0+\sigma_n \ is\ a\ cut\ time\ of\ S^1\right\},$ and
$$B_n=\left\{\left\{S^1_{\tau_0+\sigma_n+1},S^1_{\tau_0+\sigma_n+2},\cdots\right\}\cap\left\{S^1_0,S_1^1,\cdots,S^1_{\tau_0}\right\}
           =\emptyset\right\}.$$
Clearly, $A_n=B_n.$ Then for any $n\geq 1,$
\begin{eqnarray*}
&&\mathbb{P}\left(\bigcap\limits_{m=n}^{\infty}A_m\right)=\mathbb{P}\left(\bigcap\limits_{m=n}^{\infty}B_m\right)
                 =\mathbb{P}\left(B_n\right)\\
&&\ \ \ \ =\mathbb{P}\left(\left\{S_{\tau_0+\sigma_n+1}^1,S_{\tau_0+\sigma_n+2}^1,\cdots\right\}
    \cap \left\{S_0^1,S_1^1,\cdots, S_{\tau_0}^1\right\}=\emptyset\right) \\
&&\ \ \ \ \geq \mathbb{P}\left(\left\{S_{\tau_{n-1}+1}^1,S_{\tau_{n-1}+2}^1,\cdots\right\}
             \cap \left\{S_0^1,S_1^1,\cdots, S_{\tau_0}^1\right\}=\emptyset\right).
\end{eqnarray*}
Notice that SRW on $G$ is transient and $\tau_0$ is a.s. finite. Then almost surely, when $n$ is large enough,
$$\left\{S_{\tau_{n-1}+1}^1,S_{\tau_{n-1}+2}^1,\cdots\right\}
             \cap \left\{S_0^1,S_1^1,\cdots, S_{\tau_0}^1\right\}=\emptyset.$$
Therefore,
\begin{eqnarray*}
\liminf\limits_{n\rightarrow\infty}\mathbb{P}\left(\bigcap\limits_{m=n}^{\infty}A_m\right)\geq
   \mathbb{P}\left(\liminf\limits_{n\rightarrow\infty}\left\{\left\{S_{\tau_{n-1}+1}^1,S_{\tau_{n-1}+2}^1,\cdots\right\}
             \cap \left\{S_0^1,S_1^1,\cdots, S_{\tau_0}^1\right\}=\emptyset\right\}\right)=1.
\end{eqnarray*}
And further
$$\mathbb{P}(A_n\ i.o.)\geq \liminf\limits_{n\rightarrow\infty}\mathbb{P}\left(\bigcap\limits_{m=n}^{\infty}A_m\right)=1;$$
which implies that the SRW $S^1$ has a.s. infinitely many cut times.\qed
\vskip 4mm

{\bf So far we have completed proving Theorem 1.3.}\qed\\

\appendix
\section{Appendix} \label{appendix}
For readers' convenience, in this appendix, we recall a completely structural classification of quasi-transitive infinite graphs with polynomial volume growth
from \cite{WW2000}.
\vskip 4mm

\begin{defi}
Let $(X, d)$ and $\left(X^{\prime}, d^{\prime}\right)$ be two metric spaces. A rough isometry is a mapping $\varphi: X\longrightarrow X^{\prime}$ such that
for all $x,y\in X,$
\[A^{-1}d(x,y)-A^{-1}B\leq d^{\prime}(\varphi (x), \varphi (y))\leq Ad(x,y)+B, \]
and $d^{\prime}\left(x^{\prime}, \varphi (X)\right)\leq B$ for all $x^{\prime}\in X^{\prime},$ where $A\geq 1$ and $B\geq 0.$
In this case, say two spaces are roughly isometric. Particularly, if $B=0,$ then say that they are metrically equivalent.
\end{defi}
\vskip 2mm

Note that to be roughly isometric is an equivalence relation between metric spaces (\cite{WW2000} p.28 paragraph 2);
and every quasi-transitive infinite graph is roughly isometric with a vertex-transitive infinite graph (\cite{WW2000} p.29 Proposition $3.9$).
\vskip 4mm

\begin{thm}\label{Thm3}\textup{(\cite{WW2000} p.30 Theorem $3.10$)}
Assume $G$ and $G^{\prime}$ are connected infinite graphs with bounded vertex degrees, and $G$ is roughly isometric to $G^{\prime}.$ Then $G$ is recurrent $\Longleftrightarrow$ $G^{\prime}$ is recurrent. Equivalently, $G$ is transient $\Longleftrightarrow$ $G^{\prime}$ is transient.
\end{thm}
\vskip 1mm
\begin{thm}\label{Thm4}\textup{(\cite{WW2000} p.54 Theorem $5.11$)}
Let $G$ be a quasi-transitive infinite graph with volume growth function satisfying $V_G(n)\leq Cn^d$ for infinitely many $n,$
where $d$ and $C$ are two positive constants. Then $G$ is roughly isometric with a Cayley graph of some finitely generated nilpotent group.
In particular, there are a natural number $D$ and two positive constants $c_0$, $c_1$ such that
\[c_0n^D\leq V_G(n)\leq c_1(n+1)^D.\]
\end{thm}
\vskip 2mm

\begin{thm}\label{Thm6}\textup{(\cite{WW2000} p.32 Lemma $3.13$)}
Let $G$ and $G^{\prime}$ be two roughly isometric connected infinite graphs with bounded vertex degrees. Then $G$ and $G^{\prime}$ have equivalent volume growth functions, in the sense that there are two positive constants $c_0, k_0$ such that
\[V_G(n)\leq c_0V_{G^{\prime}}(k_0n)~and~V_{G^{\prime}}(n)\leq c_0V_G(k_0n)~for\ all\ n.\]
\end{thm}
\vskip 2mm

By Theorems \ref{Thm3}-\ref{Thm6}, any quasi-transitive connected infinite graph $G$ with polynomial volume growth
is roughly isometric with a Cayley graph $G^{\prime}$ of some finitely generated nilpotent group which has equivalent volume growth function with $G;$ and $G$ is recurrent (resp. transient) if and only if so is $G^{\prime}.$
Let $D$ be the degree for the polynomial volume growth of both $G$ and $G^{\prime}.$
It is known that SRW on $G^{\prime}$ is recurrent if $D\in\{1,2\},$ and transient if $D\geq 3.$
Therefore, SRW on $G$ is recurrent if $D\in\{1,2\},$ and transient if $D\geq 3.$

\end{document}